\documentclass[fleqn,3p,times]{elsarticle}
\biboptions{sort&compress}

\RequirePackage{amscd, amsfonts, amssymb}
\RequirePackage{mathrsfs, mathtools}
\RequirePackage{stmaryrd, array}
\RequirePackage{bm, graphicx, color, subfig}
\RequirePackage{url}
\RequirePackage{microtype}

\journal{Applied Mathematics Letters}

\graphicspath{{./Parts/Images/}}
\def\pathImages{./}%

\begin{document}

\newcommand{\defeq}{\vcentcolon=}

\begin{frontmatter}

\title{{Dual virtual element method in presence of an inclusion}}


\author{Alessio Fumagalli}

\address{}

\begin{abstract}
    We consider a Darcy problem  for saturated porous media written in dual
    formulation in presence of a fully immersed inclusion. The lowest order virtual element
    method is employ to derive the discrete approximation. In the present work
    we study the effect of cells with cuts on the numerical solution, able to geometrically handle in
    a more natural way the inclusion tips. The numerical
    results show the validity of the proposed approach.
\end{abstract}

\begin{keyword}
    VEM, Darcy flow, fractured porous media, inclusion
\end{keyword}

\end{frontmatter}

\section{Introduction}
Single-phase flow in fractured porous media is a challenging problem involving
different aspects, \textit{i.e.} the derivation of proper mathematical models to describe
fracture surrounding porous media flow and the subsequent
discretization with ad-hoc numerical schemes. One of the most common
approach is to consider
fractures as co-dimensional objects and derive proper reduced models and
coupling condition to describe the flow in the new setting. A hybrid
dimensional description of the problem is thus introduced, see
\cite{Karimi-Fard2004,Martin2005,DAngelo2011,Berrone2013,Formaggia2012}.
In presence of multiple fractures, forming a complex system of network,
grid creation may become challenging and the number of cells or their
shape may not be satisfactory for complex application, \textit{e.g.}, the
benchmark study proposed in \cite{Flemisch2017}.

In the present work we simplify the problem considering a single fracture
and substituting with an inclusion, \textit{i.e.} the flow problem in the fracture
is not considered but its effect of its normal flow. The inclusion is an
internal condition for the problem.  We consider two
extreme cases: perfectly permeable fracture (infinite normal permeability)
and impermeable fracture (zero normal
permeability). The pressure imposed on both sides of the inclusion determines these
two cases.

The virtual element method (VEM), introduced in
\cite{BeiraodaVeiga2013,BeiraodaVeiga2014a,Brezzi2014,BeiraodaVeiga2014b,Antonietti2014a,Benedetto2014,BeiraoVeiga2016,Benedetto2016,Fumagalli2016b},
is able to discretize the problem on grids with rather general cell shape. The theory
developed in the aforementioned works considers star shaped cells. In the
present study, the lowest order VEM is considered in presence of cells with an
internal cut for the approximation at the tips of the inclusion.

The paper is organised as follow. In Section \ref{sec:model} the
mathematical model and its weak formulation are presented. Section
\ref{sec:discrete} introduces the discrete formulation of the problem.
Numerical examples are reported in Section \ref{sec:example} for both the
extreme cases. The work finishes with conclusions Section
\ref{sec:conclusion}.

\section{Mathematical model}\label{sec:model}

\begin{figure}[htb]
    \centering
    \resizebox{0.2475\textwidth}{!}{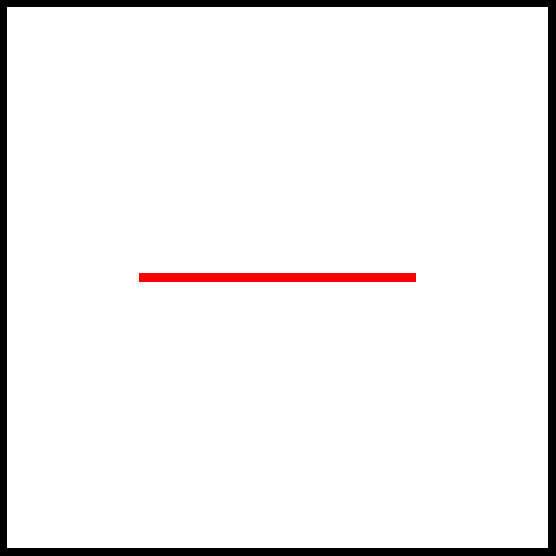}
    \caption{domain}%
    \label{fig:domain}
\end{figure}

Let us set $\Omega \subset \mathbb{R}^2$ a regular domain representing a
saturated porous media.
We consider the Darcy model for single phase flow in a saturated porous media written in dual
formulation, namely
\begin{subequations}\label{eq:darcy_inclusion}
\begin{gather}\label{eq:darcy}
    \bm{u} + K \nabla p = \bm{0} \text{ in } \Omega \quad \land \quad
    \nabla \cdot \bm{u} = f \text{ in } \Omega \quad \land \quad p = 0 \text{ on } \partial
    \Omega.
\end{gather}
The unknowns are: $\bm{u}$ the Darcy velocity and $p$ the fluid pressure.
In \eqref{eq:darcy} $K$ represents the permeability matrix, symmetric and
positive defined, and $f$ a scalar
source or sink term. To keep the presentation simple we consider only
homogeneous boundary condition for the pressure at the outer boundary of
$\Omega$, denoted by $\partial \Omega$. Coupled to \eqref{eq:darcy} we are
interested to model an immersed inclusion $\gamma$, see Figure \ref{fig:domain}
as an example. With an abuse of notation, $\partial \Omega$ does not include
$\gamma$ which represents an internal boundary for $\Omega$. It is possible to
define a unique normal $\bm{n}$ associated to $\gamma$ and two different sides
of $\gamma$ with respect to the direction of $\bm{n}$. We indicate them as $\gamma^+$
and $\gamma^-$, with $\bm{n}^+=\bm{n}$ and $\bm{n}^- = -\bm{n}$ the associated
normals. It is important to note that geometrically $\gamma$, $\gamma^+$,
and $\gamma^-$ are indeed the same object but introducing the two sides help us
to impose different internal conditions. We have
\begin{gather}\label{eq:inclusion}
    p = p^+ \text{ on } \gamma^+ \quad \land \quad
    p = p^- \text{ on } \gamma^-.
\end{gather}
\end{subequations}
We define the Hilbert spaces $Q=L^2(\Omega)$ and $V=H_{\rm
div}(\Omega)$.
By standard arguments the weak formulation of
\eqref{eq:darcy_inclusion} reads
\begin{gather}\label{eq:weak}
    a(\bm{u}, \bm{v}) + b(\bm{v}, p) = J(\bm{v}) \quad
    \forall \bm{v} \in V
    \quad \land \quad
    b(\bm{u}, q) = F(q) \quad \forall q \in Q.
\end{gather}
In \ref{eq:weak} we have indicated by
$(\cdot, \cdot)_{\Omega}$: $Q \times Q \rightarrow
\mathbb{R}$ the scalar product in $Q$. The bilinear forms and functionals are defined as
\begin{gather*}
    a(\cdot, \cdot): V \times V \rightarrow \mathbb{R}
    \quad \text{s.t.} \quad a(\bm{w}, \bm{v}) \defeq (K^{-1}\bm{w}, \bm{v})_{\Omega},
    \qquad
    b(\cdot, \cdot): V \times Q \rightarrow \mathbb{R} \quad \text{s.t.} \quad
    b(\bm{v}, q) \defeq - (\nabla \cdot \bm{v}, q)_\Omega,\\
    J(\cdot): V \rightarrow \mathbb{R} \quad \text{s.t.} \quad
    J(\bm{v}) \defeq -
    \langle \bm{v}^+ \cdot \bm{n}^+, p^+ \rangle_{\gamma^+}
    - \langle \bm{v}^- \cdot \bm{n}^-, p^- \rangle_{\gamma^-}, \quad
    F(\cdot): Q \rightarrow \mathbb{R} \quad \text{s.t.} \quad F(q) \defeq -(f, q)_{\Omega}
\end{gather*}
where the duality pairings are defined as $\langle
\cdot, \cdot \rangle_{\gamma^\pm}: H^{-\frac{1}{2}}(\gamma^\pm)
\times H^{\frac{1}{2}}(\gamma^\pm) \rightarrow \mathbb{R}$.
We are assuming  $p^\pm \in H^{\frac{1}{2}}(\gamma^\pm)$ and $f \in
L^2(\Omega)$.
Following \cite{Fumagalli2016b} problem \eqref{eq:weak} is well posed.

\section{Discrete formulation}\label{sec:discrete}

The creation of a computational grid with multiple intersecting inclusions, or
more generally fractures, is a challenging aspect possibly resulting in a high
number and/or poorly shaped cells. To overcome this aspect we consider a VEM
formulation with a clustering technique from a finer triangular grid.

We consider now the discrete formulation of problem \eqref{eq:weak} employing
the VEM. For simplicity we limit our analysis to the lowest order
case. For a more general case we refer to
\cite{BeiraodaVeiga2013,BeiraodaVeiga2014a,Brezzi2014,BeiraodaVeiga2014b,Antonietti2014a,Benedetto2014,BeiraoVeiga2016,Benedetto2016}.
Let be $\mathcal{T}(\Omega)$ a generic tessellation of $\Omega$ made of
non-overlapping polytopes. We allow non-star shaped cells by considering cells
with an internal cut. We introduce the following approximation spaces for
$p$ and $\bm{u}$. For a cell $E \in \mathcal{T}(\Omega)$ with edges
$\mathcal{E}(E)$, let us define
\begin{gather}\label{eq:discrete_spaces}
    Q_h(E)
    \defeq \{ q \in Q: \, q \in \mathbb{P}_0(E)\}
    \! \! \quad \land \quad \! \!
    V_h(E) \defeq \left\{ \bm{v} \in V: \bm{v} \cdot \bm{n} |_e \in
    \mathbb{P}_0(e)\, \forall e \in \mathcal{E}(E), \, \nabla \cdot \bm{v} \in
    \mathbb{P}_0(E), \nabla \times \bm{v} = \bm{0}
    \right\}.
\end{gather}
The norms and the global spaces, $Q_h$ and $V_h$, are defined accordingly. With
the definition \eqref{eq:discrete_spaces} it is possible to compute immediately
the discrete approximation of $b(\cdot, \cdot)$, $J(\cdot)$, and $F(\cdot)$.
Since the shape of the functions in $V_h(E)$, for each cell $E \in
\mathcal{T}(E)$, is not prescribed we cannot compute directly the bilinear form
$a(\cdot, \cdot)$. In this case we need to introduce a sub-space
of $V_h$ and a projection operator, the space is
\begin{gather*}
    \mathcal{V}_h(E) \defeq \left\{ \bm{v} \in V_h(E): \bm{v} = \nabla v, \text{
    for } v \in \mathbb{P}_1(E)\right\}
\end{gather*}
It is worth noting that the local approximation of $\mathbb{P}_1(E)$ is done by
means of a suitable monomial expansion. The projection operator is defined locally as
$\Pi_0 : V(E) \rightarrow \mathcal{V}_h(E)$ such that $a(T_0 \bm{u}, \bm{v}) =
0$ for all $\bm{v} \in \mathcal{V}_h(E)$, with $T_0 \defeq I - \Pi_0$ the
projection on the orthogonal space of $\mathcal{V}_h$. Considering the property
of the projection operator, we make the following approximation
$a(\bm{u}, \bm{v}) \approx a_h(\bm{u}, \bm{v}) \defeq a(\Pi_0 \bm{u}, \Pi_0
\bm{v} ) + s(T_0 \bm{u}, T_0 \bm{v})$, where the fist part is a consistency term
and the second a stabilization term. The latter is approximated by the bilinear
form $s(\cdot, \cdot)$ such that an equivalence property holds true,
\textit{i.e.}
$\exists \iota_*, \iota^* \in \mathbb{R}$ independent from the
discretization size $h$ satisfying $\iota_* a(\Pi_0 \bm{v}, \Pi_0 \bm{v}) \leq
s(T_0 \bm{v}, T_0 \bm{v}) \leq \iota^* a(\Pi_0 \bm{v}, \Pi_0 \bm{v})$,
see the aforementioned work for more details.
With these choices the bilinear form $a_h(\cdot, \cdot)$ is computable for
each cell of the grid.



\section{Numerical example}\label{sec:example}

In this section we present a numerical study to investigate the error decay for
the model presented previously. Two examples are shown with
the same and different boundary conditions on the internal boundary $\gamma$.
We consider the domain depicted in
Figure \ref{fig:domain} with $\Omega=[0, 1]^2$ and $\gamma=\{(x, y): 0.25 \leq x
\leq 0.75 \land y =0.5\}$.

The reference solution, named $p_{\rm ref}$, is computed considering a grid of
77546 triangles very refined at the ending points of $\gamma$. The $L^2$ relative
error presented below is defined as
\begin{align}\label{eq:error}
    err(p) = \dfrac{||p - p_{\rm ref}||_{L^2(\Omega)}}{|\max{p_{\rm ref}} -
    \min{p_{\rm ref}}|}.
\end{align}

\begin{figure}[tb]
    \centering
    \includegraphics[width=0.2475\textwidth]{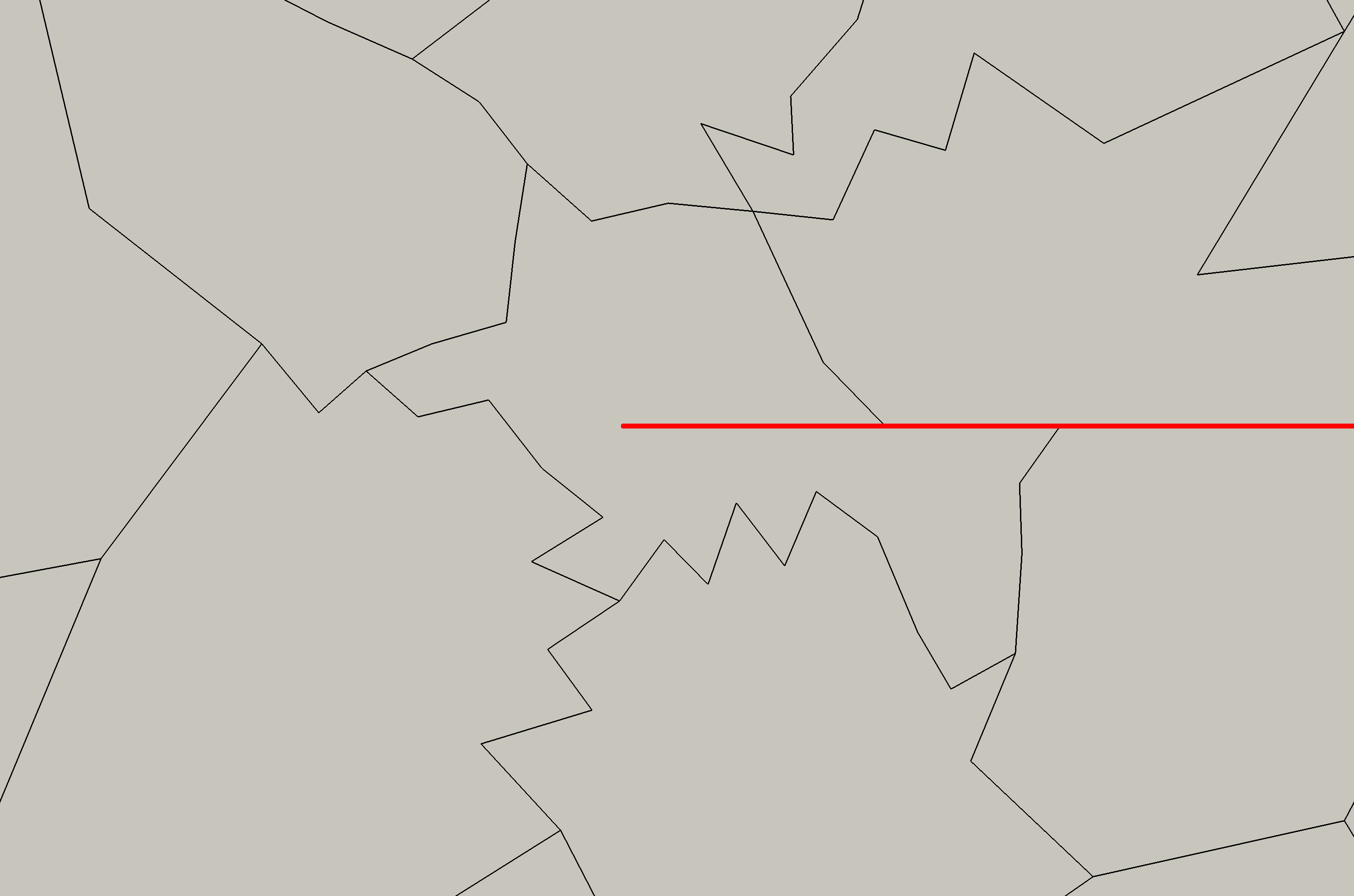}\hfill%
    \includegraphics[width=0.2475\textwidth]{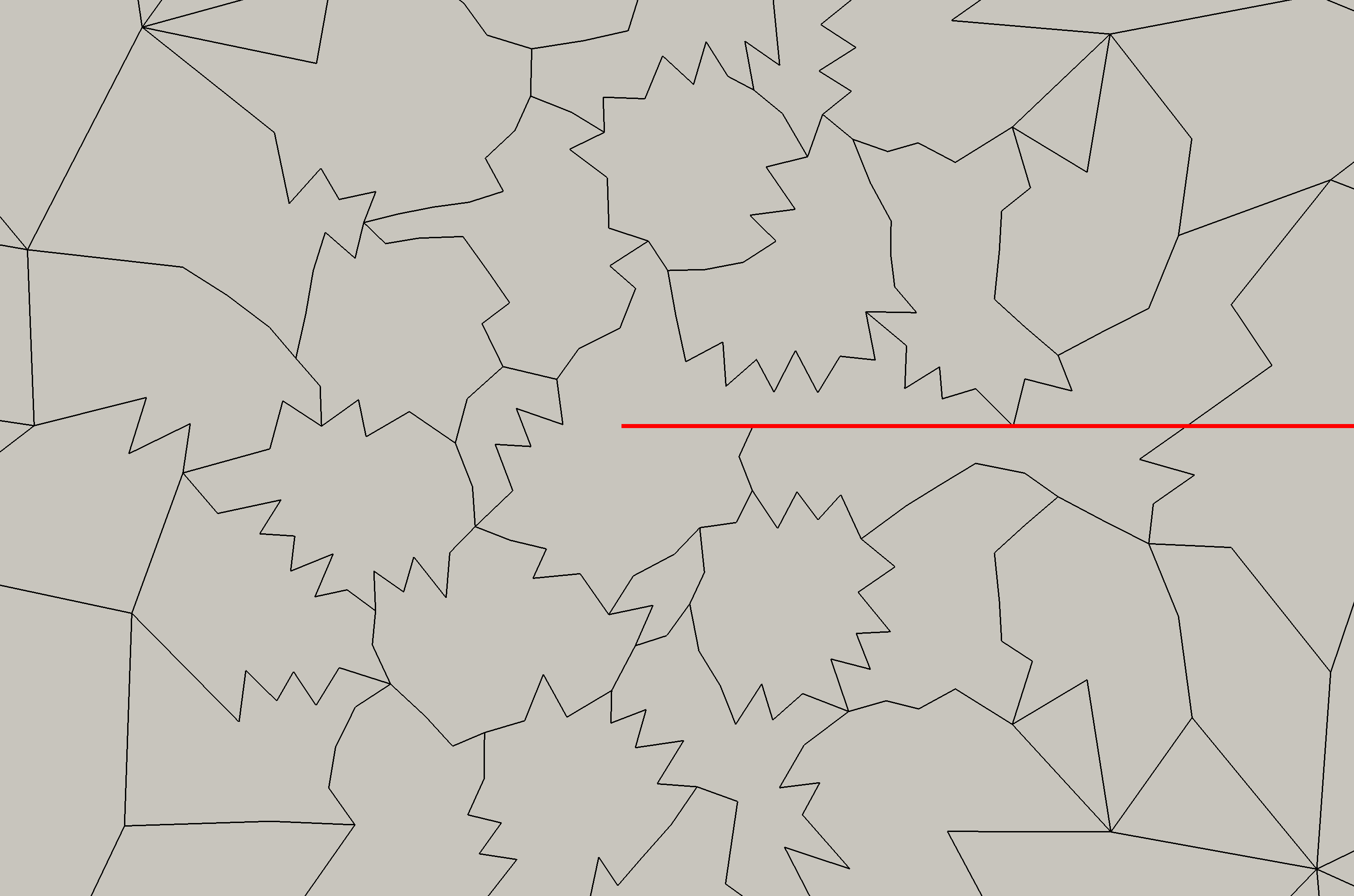}\hfill%
    \includegraphics[width=0.2475\textwidth]{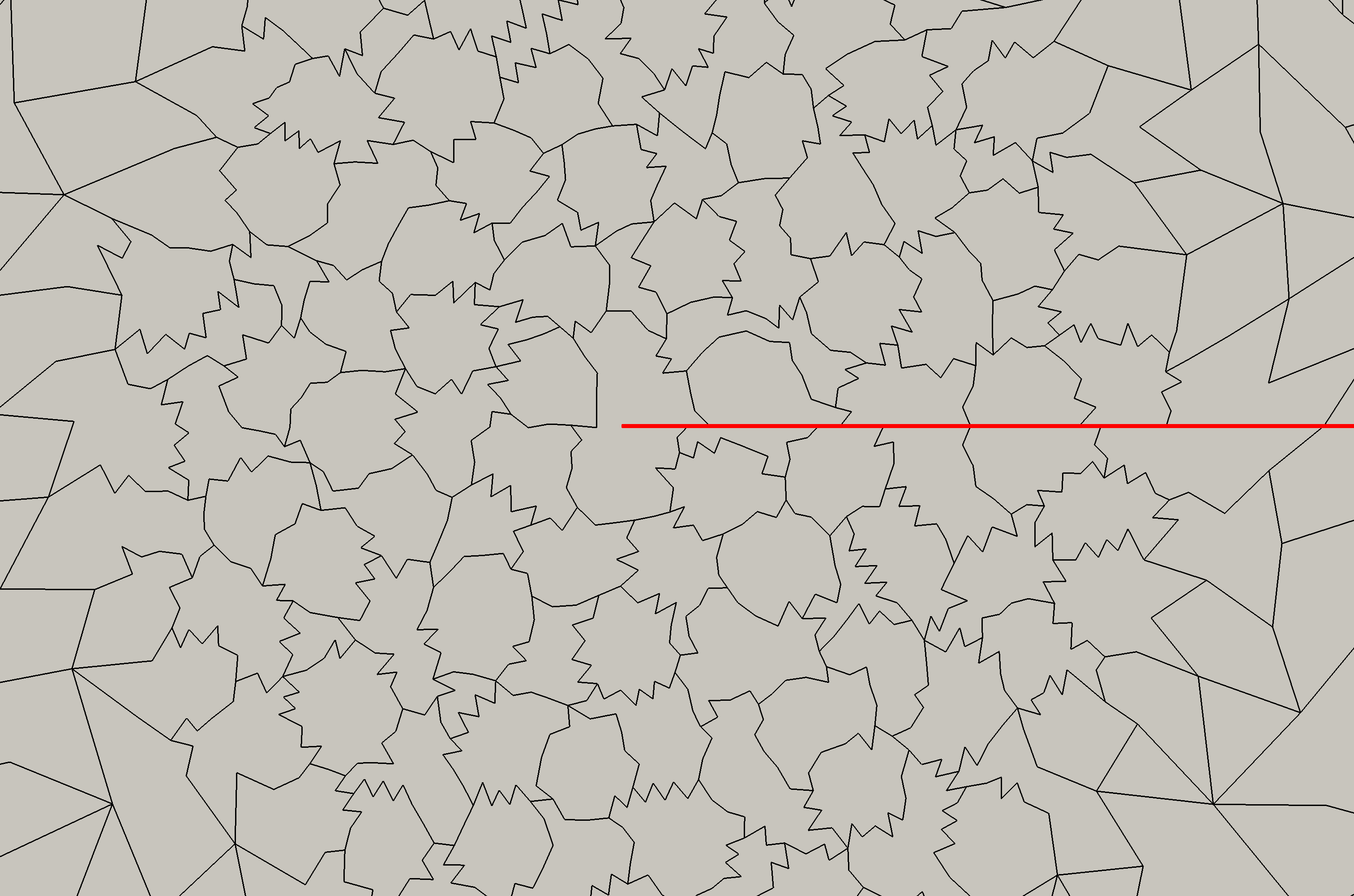}\hfill%
    \includegraphics[width=0.2475\textwidth]{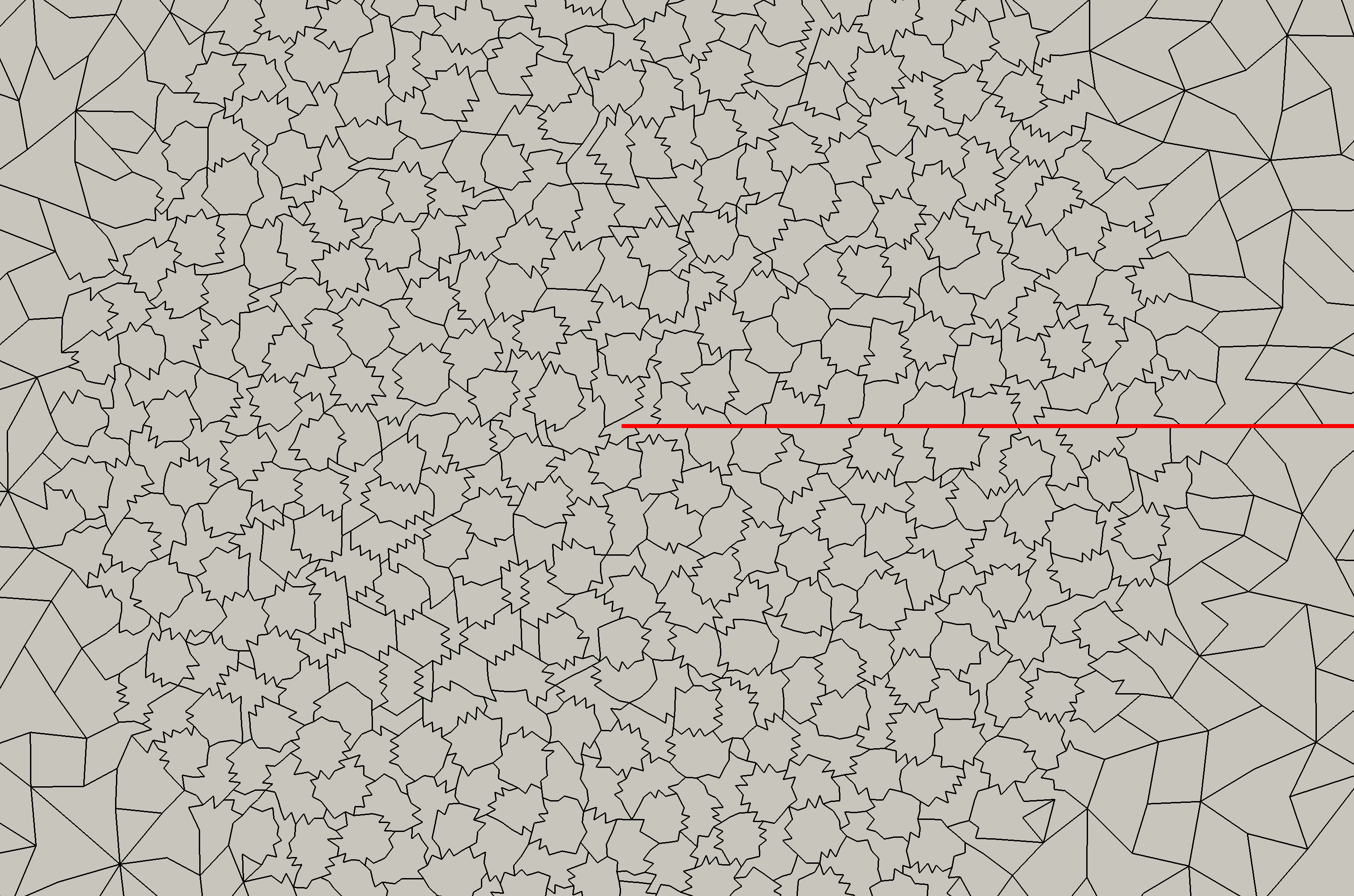}
    \caption{Detail of the computational grids on the left tip of
        $\gamma$, depicted in red.}%
    \label{fig:zoom}
\end{figure}

In both cases we generate a family of triangular grids with a different level of
refinement, in particular, close to the ending points of $\gamma$. Later, an
agglomeration technique is employed based on the cell measure: neighbour cells
with small measure are glue together forming new cells. To stress the method
presented previously we enforce the creation, at the tips of $\gamma$, of cells
with internal cuts. Figure \ref{fig:zoom} shows a zoom of the considered grids
at the left ending point of $\gamma$.  It is worth noting that the computation
of the mesh size $h$ is itself complex.

The implementation is done with PorePy: a simulation tool for fractured
and deformable porous media written in Python. See
\url{github.com/pmgbergen/porepy} for further information.
All the examples of this section are available in the package.
In the following pictures a ``Blue to Red Rainbow'' colour map is used.


\subsection{Continuous pressure condition} \label{subsec:cont}

\begin{figure}[tb]
    \centering
    \subfloat[Reference and coarse solution for the continuous case.]
    {
        \centering
        \includegraphics[width=0.2475\textwidth]{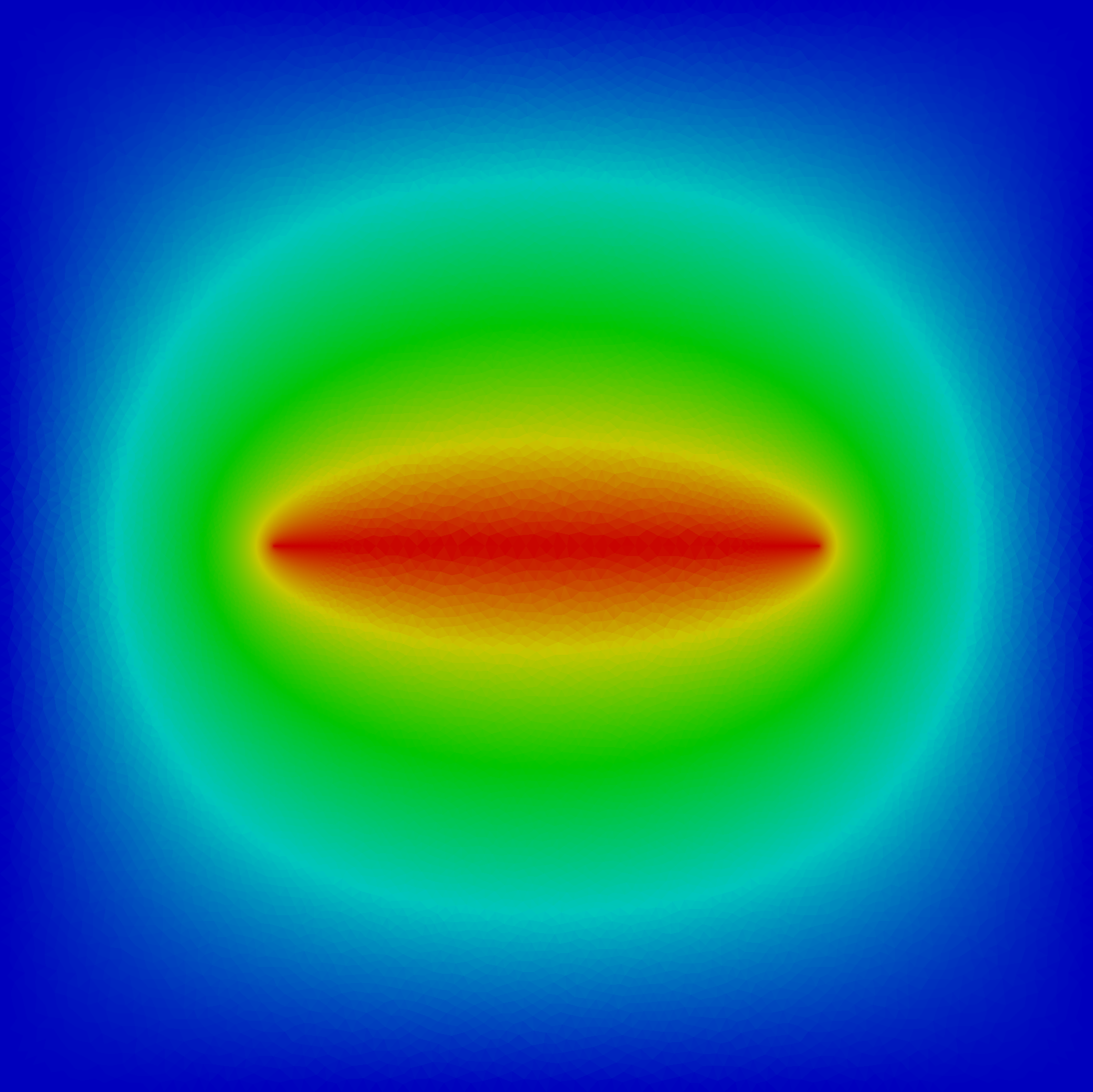}%
        \includegraphics[width=0.2475\textwidth]{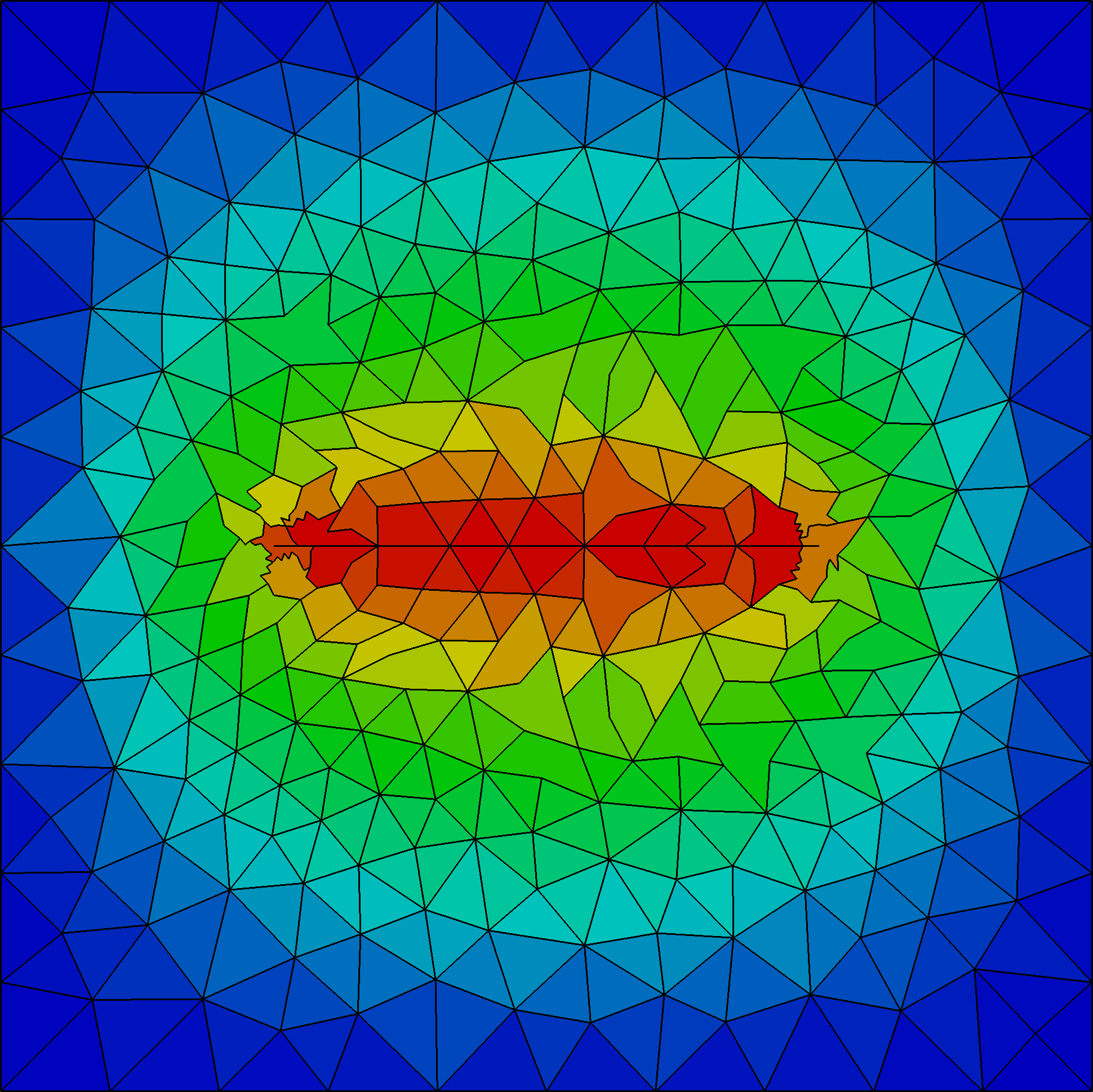}%
    }
    \subfloat[Reference and coarse solution for the discontinuous case.]
    {
        \centering
        \includegraphics[width=0.2475\textwidth]{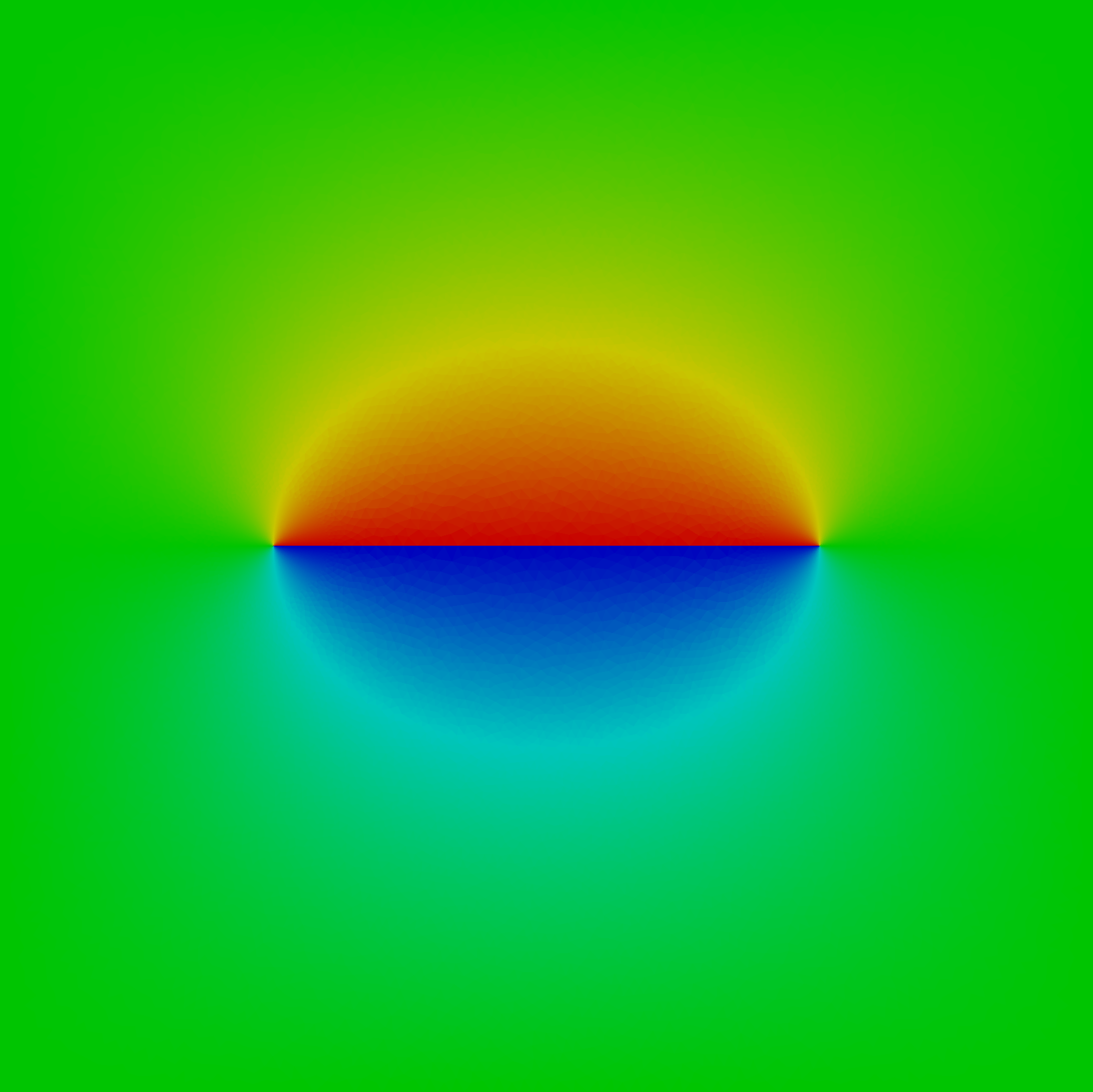}%
        \includegraphics[width=0.2475\textwidth]{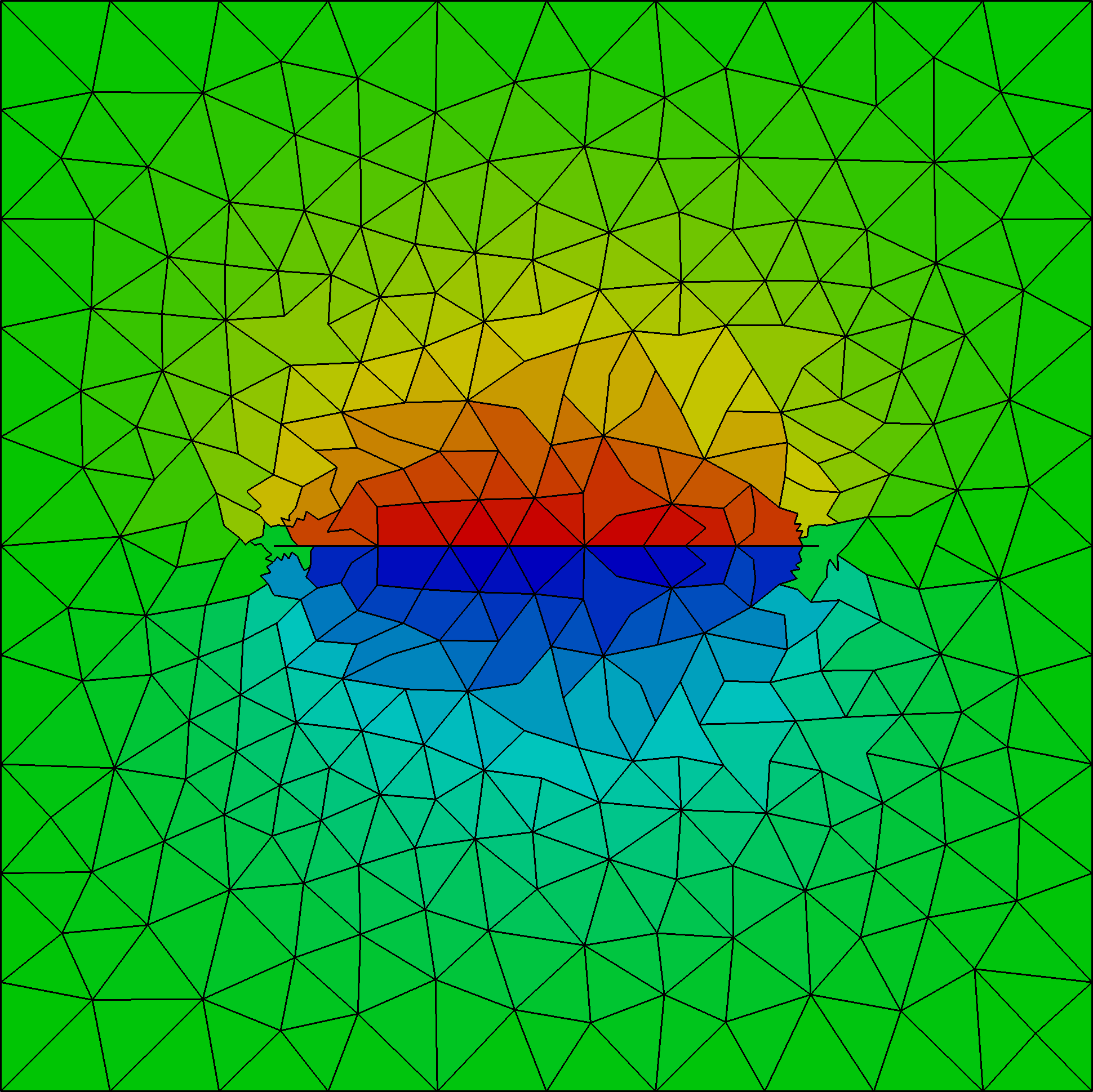}%
    }%
    \caption{On the left: pressure (range $(0, 1)$) for the
        continuous case. On the right: pressure (range $(-1, 1)$) for
        the discontinuous case.}%
    \label{fig:continuous}
\end{figure}

\begin{figure}[tb]
    \centering
    \includegraphics[width=0.3475\textwidth]{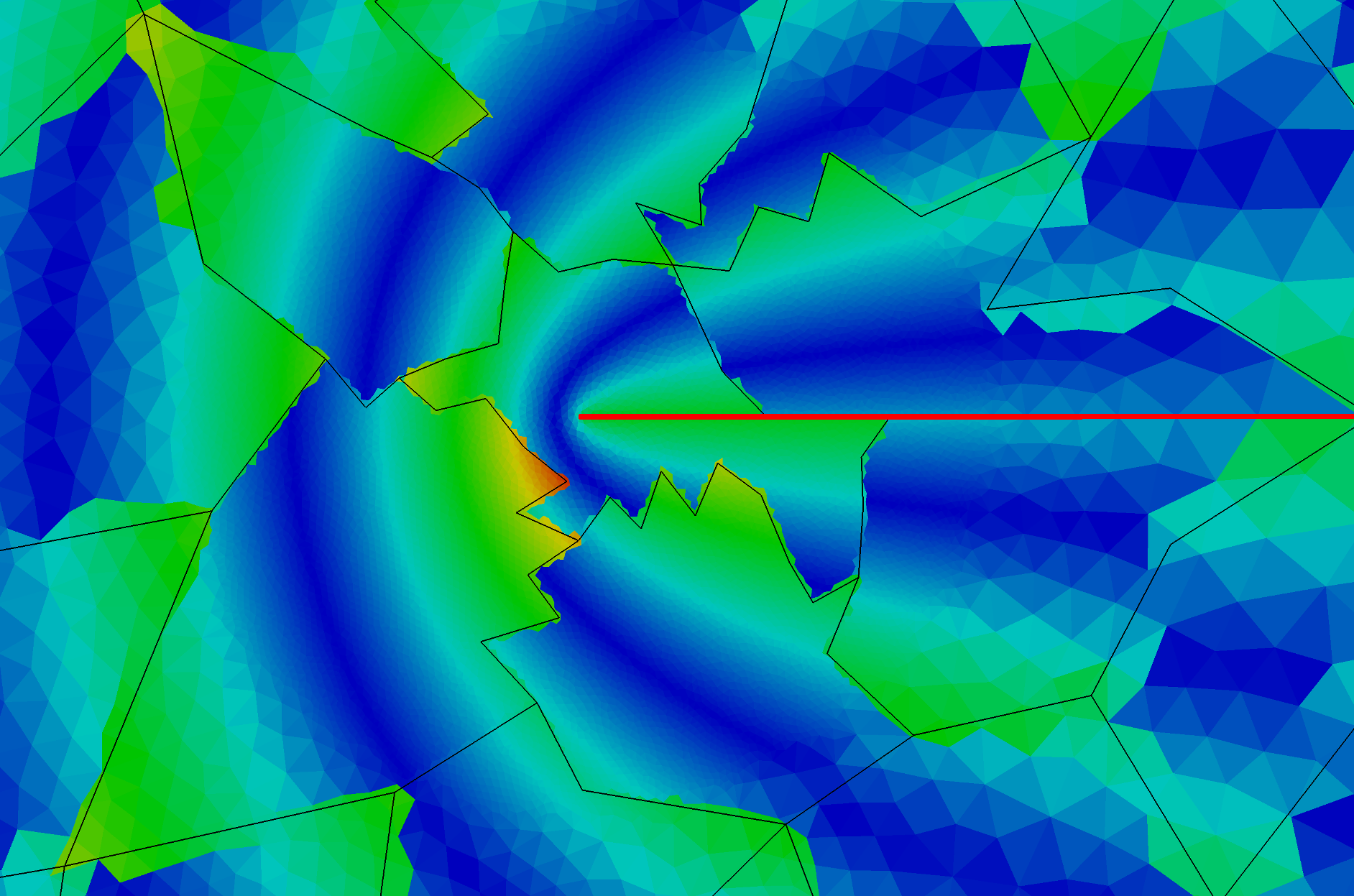}%
    \hspace*{0.025\textwidth}%
    \includegraphics[width=0.25\textwidth]{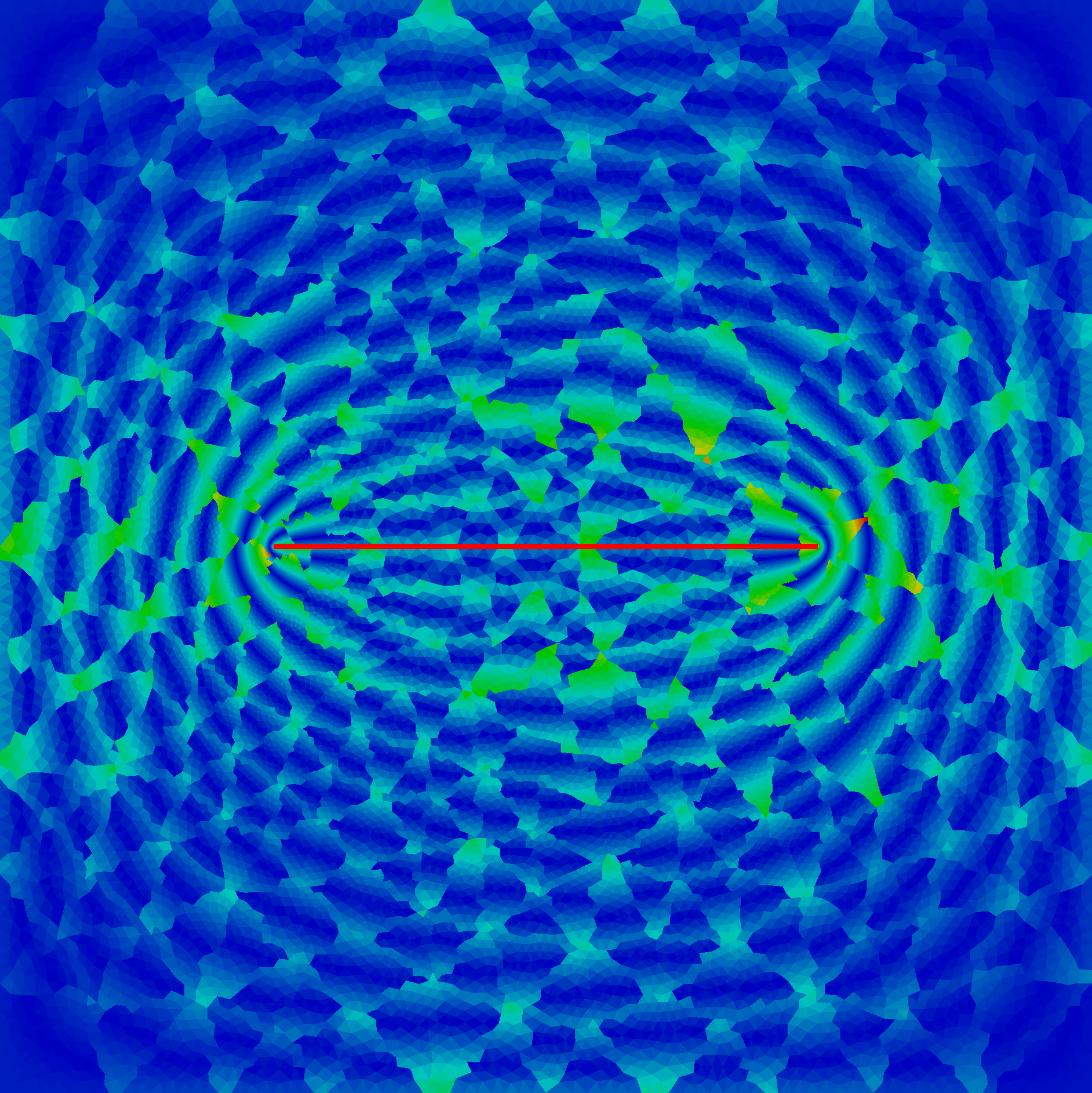}%
    \hspace*{0.025\textwidth}%
    \includegraphics[width=0.3475\textwidth]{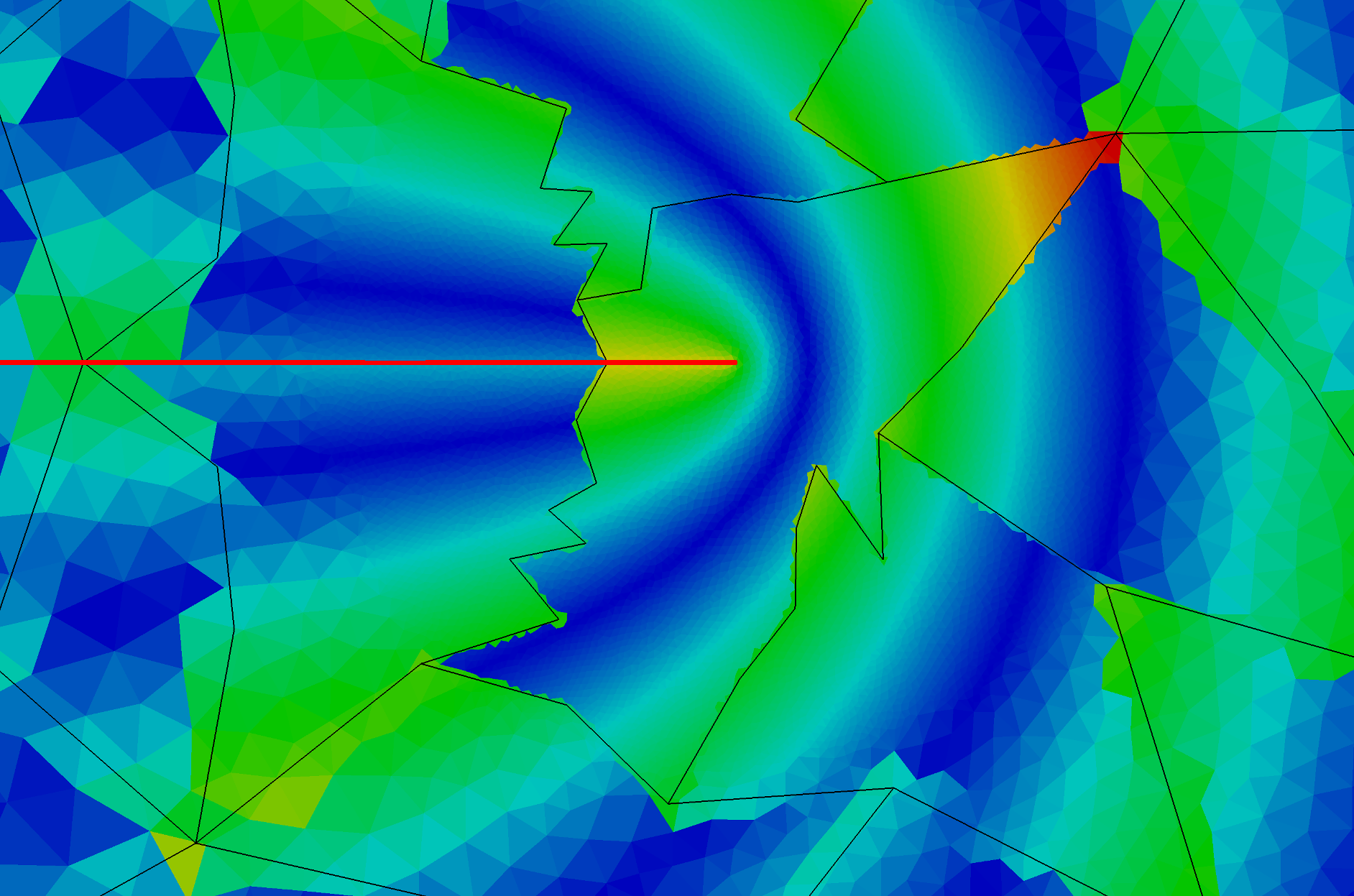}
    \caption{In the centre: the error on $\Omega$ represented on the reference
        grid. On the right and left: a zoom of the error for the right and left tip,
        respectively. Range in $(0, 0.25)$.}%
    \label{fig:continuous_error}
\end{figure}

We consider $p^+ = p^- = 1$ as condition for $\gamma^+$ and $\gamma^-$.
The pressure solution of the reference and of the coarser grid is reported in
Figure \ref{fig:continuous} on the left.
We notice that the solution computed on
the coarse grid with cells containing a cut resembles the reference solution,
even at the tips of the inclusion. The error, computed for each cell of the
reference grid, is represented in Figure \ref{fig:continuous_error}.
The error
is distributed in the domain without any specific peak at the tips. In this case
we can conclude that the cells with a cut do not deteriorate the quality of the
solution. Finally, the error decay is given in Table \ref{tab:error} on the
left confirming a first order of convergence.
\begin{table}
    \centering
    \begin{tabular}{|c|c|c||c|c|}
        \hline
         $h$ & $err(p)$ & $\mathcal{O}(h)$ & $err(p)$ & $\mathcal{O}(h)$\\ \hline
         7.942$\cdot 10^{-2}$ & 3.045$\cdot 10^{-2}$ & -    & 2.266$\cdot 10^{-2}$ & -      \\ \hline
         3.801$\cdot 10^{-2}$ & 1.545$\cdot 10^{-2}$ & 0.92 & 1.217$\cdot 10^{-2}$ & 0.84   \\ \hline
         1.920$\cdot 10^{-2}$ & 8.03 $\cdot 10^{-3}$ & 0.96 & 5.916$\cdot 10^{-3}$ & 1.06   \\ \hline
         9.692$\cdot 10^{-3}$ & 1.830$\cdot 10^{-3}$ & 2.16 & 2.913$\cdot 10^{-3}$ & 1.04   \\ \hline
    \end{tabular}
    \caption{On the left: \eqref{eq:error} for example in
        Subsection \ref{subsec:cont}.
        On the right: \eqref{eq:error} for the example in Subsection
        \ref{subsec:discont}.}%
    \label{tab:error}
\end{table}


\subsection{Discontinuous pressure condition} \label{subsec:discont}

\begin{figure}[tb]
    \centering
    \includegraphics[width=0.3475\textwidth]{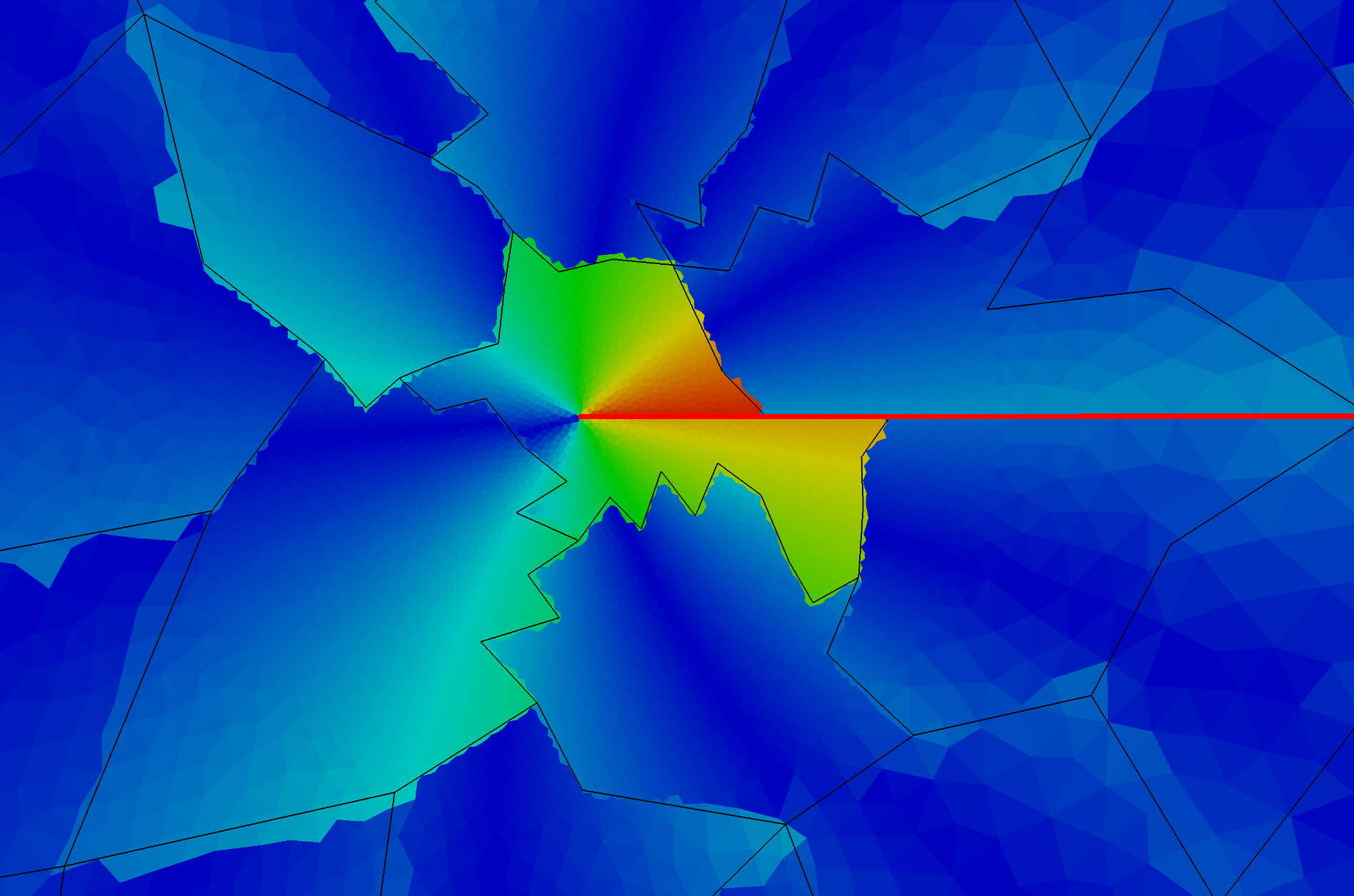}%
    \hspace*{0.025\textwidth}%
    \includegraphics[width=0.25\textwidth]{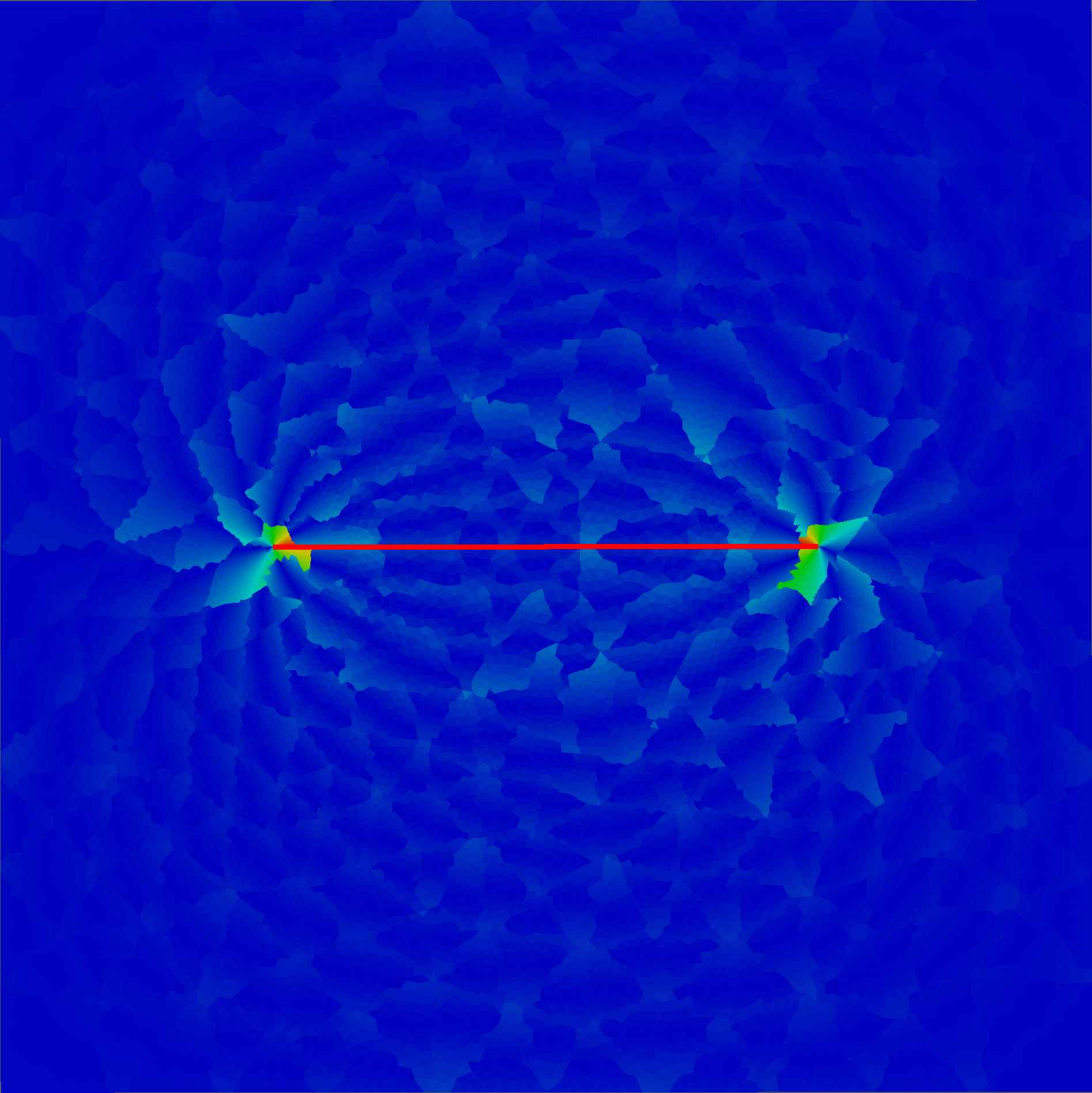}%
    \hspace*{0.025\textwidth}%
    \includegraphics[width=0.3475\textwidth]{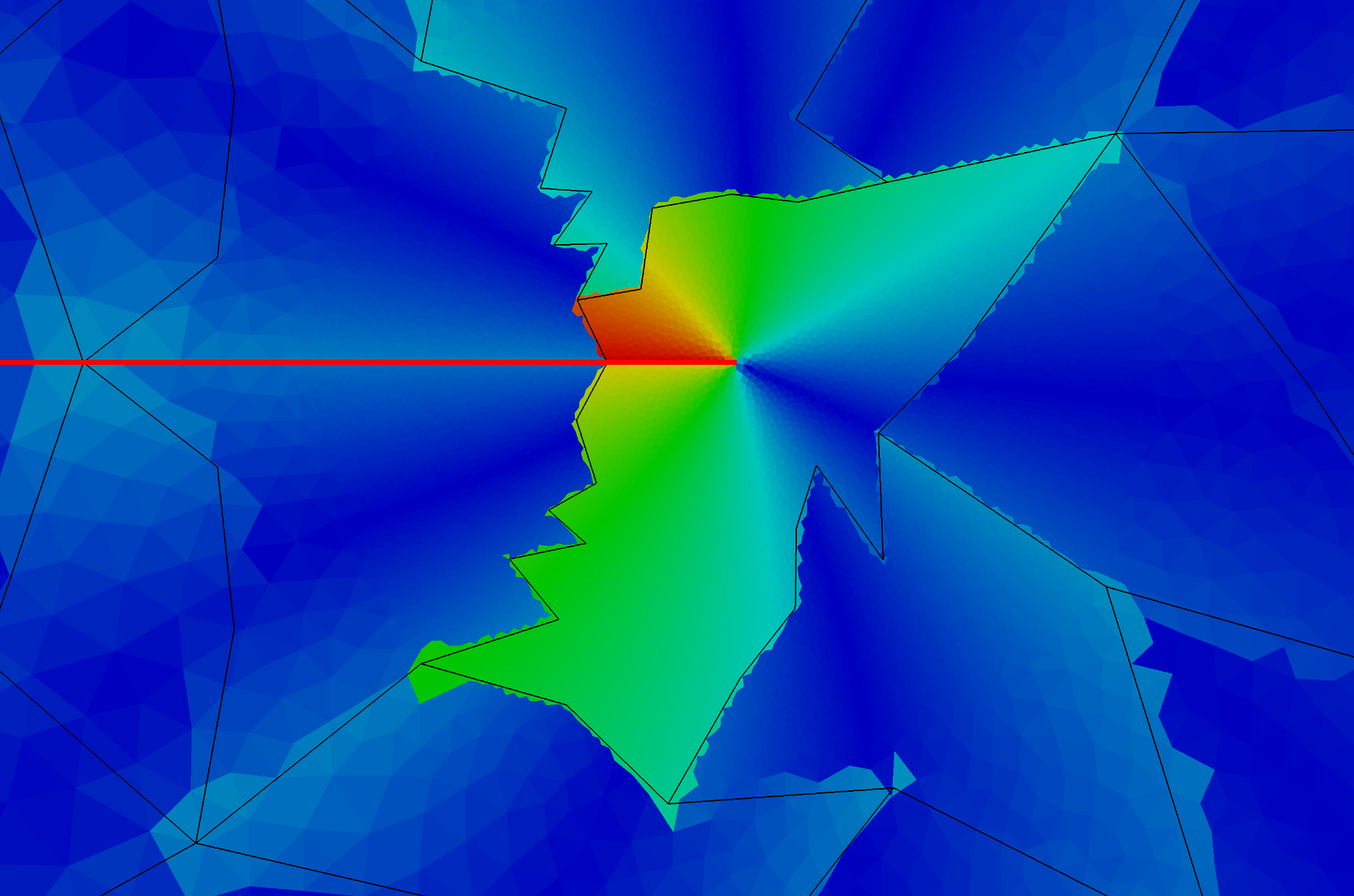}
    \caption{From the left: reference solution, coarse solution, and
    solution (range $(0, 0.57)$) for the
        continuous case.}%
    \label{fig:discontinuous_error}
\end{figure}

We consider $p^+ = 1$ for $\gamma^+$ and $p^- = -1$ for $\gamma^-$ as boundary
condition for the inclusion. We expect a jump in the solution and a poor
regularity at the tips of $\gamma$. The pressure solution of the reference and
of a coarse grid is reported in Figure \ref{fig:continuous} on the right. The
solution computed with the coarse grid is in accordance with the reference
solution, even at the tips of the inclusion. The error, computed for each cell
of the reference grid, is depicted in Figure \ref{fig:discontinuous_error}.
We point out that the error is focused at the tips of the inclusion, since a
complex solution is now approximated with a single degree of freedom on each
coarse cell containing the tip. Nevertheless, the errors listed in Table
\ref{tab:error} confirm also in this case a first order of convergence, thus the
quality of the computed solution.



\section{Conclusion}\label{sec:conclusion}

In this work we studied the applicability of a dual virtual element method in
presence of cells with a cut for a problem with an inclusion.
The proposed scheme behaved accurately obtaining correct results and error
decay. For a different conditions at the inclusion a pick of error was
concentrated at the inclusion tips, because a complex solution was approximated
with a constant. However, the solution was not deteriorate and the global error
was acceptable. While in the case of equal condition at the inclusion, a more
uniform error was observed inside the domain. This is a promising approach to
enlight the computational cost in presence of multiple inclusions, or even
fractures represented as co-dimensional domains, which is currently under
investigation. 



\section*{Acknowledgment}

    We acknowledge financial support for the ANIGMA project from the Research
    Council of Norway (project no. 244129/E20) through the ENERGIX program.
    The author wish to thank:
    Runar Berge,
    Inga Berre,
    Wietse Boon,
    Eirik Keilegavlen, and
    Ivar Stefansson for many fruitful discussions.





\section*{Bibliography}





\end{document}